\documentclass[preprint,12pt]{elsarticle}
\usepackage[nodots,nocompress]{numcompress}
\usepackage{epsfig,amsmath,amssymb,latexsym,amsfonts,curves,eepic,epsfig,epic,graphics,amsthm,graphicx}
\usepackage{hyperref}
\numberwithin{equation}{section}

\newtheorem{remark}{Remark}[section]

\usepackage{mathrsfs}
\usepackage{mathtools}
\usepackage{graphicx}
\usepackage{stmaryrd}
 \usepackage{subfigure}
\usepackage{float}
\usepackage{amsmath}
    \usepackage{bm}
\usepackage{amsfonts,amssymb}
   \usepackage{natbib}
  \usepackage{dsfont}
\usepackage{amsfonts}
\usepackage{mathrsfs}
  \usepackage{pifont}

\begin{document}
\begin{frontmatter}
\title{ {\Large
Discrete line integral method for the Lorentz force system}}
\author{ Haochen Li$^{a}$, \quad Yushun Wang$^{a,*}$ }
\address{
$^{a}$ Jiangsu Key Laboratory for NSLSCS, School of Mathematical Sciences, Nanjing Normal University, Jiangsu, $210023$, People's Republic of China\\}
\begin{abstract}
In this paper, we apply the Boole discrete line integral to solve the Lorentz force system which is written as a non-canonical Hamiltonian system.
The method is exactly energy-conserving for polynomial Hamiltonians of degree $\nu \leq 4$.
In any other case, the energy can also be conserved approximatively.
With comparison to well-used Boris method, numerical experiments are presented to demonstrate the
energy-preserving property of the method.

\end{abstract}
\end{frontmatter}

\begin{figure}[b]
\small \baselineskip=10pt
\rule[2mm]{1.8cm}{0.2mm} \par
$^{*}$Corresponding author.\\
E-mail address: wangyushun@njnu.edu.cn(Y. Wang), Tel.: +86-13912961556, Fax: +86-517-83525012.
\end{figure}

\pagestyle{myheadings}
\markboth{\hfil 
   \hfil \hbox{}}
{\hbox{} \hfil 
Haochen Li,  Yushun Wang  \hfil}

\section{Introduction}
\label{intro}

Geometric numerical integration methods have come to the fore,
partly as an alternative to traditional methods such as Runge-Kutta methods.
A numerical method is called geometric if it preserves one or more physical/geometric
properties of the system exactly (i.e. up to round-off error)~\cite{Ref1,Ref2,Ref3,Ref4}. Examples of such geometric
properties that can be preserved are first integrals, symplectic structures, symmetries and reversing symmetries, phase space volumes, Lyapunov functions, foliations, etc.
Geometric methods have applications in many areas of physics, including celestial mechanics,
particle accelerators, molecular dynamics, fluid dynamics, pattern formation, plasma physics, reaction-diffusion equations, and meteorology~\cite{Refhnw,Refhlw,Ref5,Ref6,Refkhz,Refxsc}.

Many important phenomena in plasmas can be understood and analyzed in terms of the single-particle motion which satisfies the Lorentz force equations~\cite{RefBellan}.
The motion of charged particles in single particle model is governed by the Newton equation under the Lorentz force exerted by a given electromagnetic field. Hamiltonian formulation is available for the single particle model and other magnetized plasma models in special coordinates~\cite{RefLittlejohn,Refmg}.
In long-term simulating the motion of charged particles, non-geometric methods such as the standard 4th order Runge-Kutta method may rise to a complete wrong solution orbit, since the numerical errors of each time step will add up coherently and become significantly large over many time steps.
In contrast, geometric methods show a good long-term accuracy, for instance, the volume-preserving algorithms~\cite{Refhslq}, the Boris method which can also conserve phase space volume~\cite{RefBoris,Refscpw,Refbl,Refqzxlst}, the variational symplectic method~\cite{Refqgt,Refqg,RefLittlejohn2} and so on. It is well-know that the most noticeable structure of a Hamiltonian system is the
Hamiltonian function itself which is usually the energy of the system.
The paper is devoted to constructing an energy-preserving method for the system.

The conservation of the energy function is one of the most relevant features characterizing a Hamiltonian system.
Methods that exactly preserve energy have been considered since several decades.
Many energy-preserving methods have been proposed~\cite{Ref31,Ref32,Ref37,Ref30,Ref36}.
The discrete gradient method is among the most popular methods for designing integral
preserving schemes for ordinary differential equations, which was perhaps first discussed
by O. Gonzalez~\cite{Ref31}. T. Matsuo proposed discrete variational method for nonlinear wave equation~\cite{Ref32}.
The averaged vector field method which is a B-series method has been proposed~\cite{Ref37,Ref30}.
More recently, L. Brugnano and F. Iavernaro proposed the discrete line integral (DLI) methods~\cite{Refip1,Refbi} and the Hamiltonian boundary value methods~\cite{Ref33,Ref34}.
The key idea of the DLI methods is to exploit the relation between the method itself and the discrete
line integral, i.e., the discrete counterpart of the line integral in conservative vector fields.
This tool yields exact conservation for polynomial Hamiltonians of arbitrarily high-degree.
Different qudarature formulas yield different DLI methods.
Especially, if we choose the Boole's rule which is the Newton-Cotes formula of degree 4, the so-called
Boole discrete line integral (BDLI) method which is exactly energy-conserving for a polynomial Hamiltonian of degree $\nu \leq 4$ is obtained.

In this paper, the Lorentz force system is written as a non-canonical Hamiltonian form.
We apply the BDLI method for the Hamiltonian system, and a new energy-preserving method is obatined.
The new method is symmetric and can preserve the Hamiltonian up to round-off error.

RefBoris

The paper is organized as follows: In Sect.~\ref{sec:3}, the dynamics of charged particles in the electromagnetic field is shown and it is written as a non-canonical Hamiltonian system.
In Sect.~\ref{sec:4}, we use the BDLI method to solve the Hamiltonian system. Based on the Boole's rule, a new method for Lorentz force system is obtained.
With comparison to well-known Boris method~\cite{RefBoris}, numerical experiments are
presented in Sect.~\ref{sec:5} to confirm the theoretical results.
We finish the paper with conclusions in Sect.~\ref{sec:6}.

\section{Hamiltonian form of the Lorentz force system}
\label{sec:3}

In this section, we review the Hamiltonian form of the Lorentz force system~\cite{RefLittlejohn,Refmg,Refhslq}.
For a charged particle in the electromagnetic field, its dynamics is governed by the Newton-Lorentz equation
\begin{equation}\label{eq:3s1}
m\textbf{\"{x}}=q(\textbf{E}+\textbf{\.{x}}\times \textbf{B}),\quad x \in \mathbb{R}^{3},
\end{equation}
where $\textbf{x}$ is the position of the charged particle, $m$ is the mass, and $q$ denotes the electric charge. For convenience, here we assume that $\textbf{B}$ and $\textbf{E}$ are static, thus $\textbf{B}=\nabla\times \textbf{A}$   and $\textbf{E}=-\nabla \varphi$ with $\textbf{A}$ and $\varphi$ the potentials.

Let the conjugate momentum be $\textbf{p}=m\textbf{\.{x}}+q\textbf{A}(\textbf{x})$, then the system \eqref{eq:3s1} is Hamiltonian with
\begin{equation}\label{eq:3s2}
H(\textbf{x},\textbf{p})=\frac{1}{2m}(\textbf{p}-q\textbf{A}(\textbf{x}))\cdot(\textbf{p}-q\textbf{A}(\textbf{x}))+q\varphi(\textbf{x}).
\end{equation}

Recasting \eqref{eq:3s1} with transformation $G:(\textbf{x},\textbf{p})\longrightarrow (\textbf{x},\textbf{v}), \textbf{x} =\textbf{x}, \textbf{v} =\textbf{p}/m-q\textbf{A}(\textbf{x})/m$, we can obtain
\begin{align}\label{eq:3s3}
&\textbf{\.{x}}=\textbf{v},\\ \label{eq:3s4}
&\textbf{\.{v}}=\frac{q}{m}(\textbf{E}(\textbf{x})+\textbf{v}\times \textbf{B}(\textbf{x})).
\end{align}
Denote $\textbf{z} =[\textbf{x}^{T}, \textbf{v}^{T}]^{T}$.
\eqref{eq:3s3}-\eqref{eq:3s4} can be written as a non-canonical Hamiltonian system
\begin{equation}\label{eq:3s5}
\textbf{\.{z}}=f(\textbf{z})=K(\textbf{z})\nabla H(\textbf{z}),
\end{equation}
where  $H(\textbf{z}) =m\textbf{v}\cdot \textbf{v}/2 +q\varphi(\textbf{x})$,
and
\[
K(\textbf{z})=
\left(
       \begin{array}{cc}
               0  & \frac{1}{m}I  \\
              -\frac{1}{m}I  & \frac{q}{m^{2}}\textbf{\^{B}}(\textbf{x})
             \end{array}
\right)
\]
is a skew-symmetric matrix with
\[
\textbf{\^{B}}(\textbf{x})=
\left(
       \begin{array}{ccc}
              0   & B_{3}(\textbf{x}) & -B_{2}(\textbf{x}) \\
              -B_{3}(\textbf{x})  & 0 & B_{1}(\textbf{x}) \\
              B_{2}(\textbf{x}) & -B_{1}(\textbf{x}) & 0
             \end{array}
\right),
\]
defined by $\textbf{B}(\textbf{x})=[B_{1}(\textbf{x}), B_{2}(\textbf{x}), B_{3}(\textbf{x})]^{T}$.

\section{Boole discrete line integral method}
\label{sec:4}

It is well know that the flow of the system \eqref{eq:3s5} preserves the energy which is usually the Hamiltonian $H(\textbf{z})$ exactly.
In this section, we derive a new energy-preserving scheme for the system \eqref{eq:3s5} by using the BDLI method proposed in~\cite{Refip1,Refbi}.
Starting
from the initial condition $\textbf{z}_{0}$ we want to produce a new approximation at $t = h$, say $\textbf{z}_{1}$,
such that the Hamiltonian is conserved. By considering the simplest possible path joining
$\textbf{z}_{0}$ and $\textbf{z}_{1}$, i.e., the segment
\begin{equation}\label{eq:4s1}
\sigma(ch)=c\textbf{z}_{1}+(1-c)\textbf{z}_{0},\quad c\in [0,1],
\end{equation}
we obtain that
\begin{align}\label{eq:4s2}
\frac{1}{h}(H(\textbf{z}_{1})-H(\textbf{z}_{0})) & = \frac{1}{h}(H(\sigma(h))-H(\sigma(0))) \\
& =  \frac{1}{h}\int_{0}^{h}\nabla H(\sigma(t))^{T}\sigma'(t)dt \nonumber\\
& =  \int_{0}^{1}\nabla H(\sigma(ch))^{T}\sigma'(ch)dc \nonumber\\
& =  \frac{1}{h}\int_{0}^{1}\nabla H(c\textbf{z}_{1}+(1-c)\textbf{z}_{0})^{T}(\textbf{z}_{1}-\textbf{z}_{0})dc \nonumber\\
& =  [\int_{0}^{1}\nabla H(c\textbf{z}_{1}+(1-c)\textbf{z}_{0})dc]^{T}\frac{\textbf{z}_{1}-\textbf{z}_{0}}{h}=0, \nonumber
\end{align}
provided that
\begin{equation}\label{eq:4s3}
\frac{\textbf{z}_{1}-\textbf{z}_{0}}{h}=K(\textbf{\^{z}}) \int_{0}^{1}\nabla H(c \textbf{z}_{1}+(1-c)\textbf{z}_{0})dc,
\end{equation}
where we choose $\textbf{\^{z}}=\frac{\textbf{z}_{1}+\textbf{z}_{0}}{2}$ in this paper.
In fact, due to the fact that $K(\frac{\textbf{z}_{1}+\textbf{z}_{0}}{2})$ is skew symmetric, we have
\begin{align}\label{eq:4s4}
&[\int_{0}^{1}\nabla H(c\textbf{z}_{1}+(1-c)\textbf{z}_{0})dc]^{T}\frac{\textbf{z}_{1}-\textbf{z}_{0}}{h}\\
&=[\int_{0}^{1}\nabla H(c\textbf{z}_{1}+(1-c)\textbf{z}_{0})dc]^{T}K(\frac{\textbf{z}_{1}+\textbf{z}_{0}}{2})[\int_{0}^{1}\nabla H(c\textbf{z}_{1}+(1-c)\textbf{z}_{0})dc]=0.\nonumber
\end{align}

In this paper, we use the Boole's rule to calculate the integrand at the right-hand side in \eqref{eq:4s3} numerically.
Then we obtain the BDLI method for \eqref{eq:3s5}:
\begin{align}\label{eq:4s5}
\textbf{z}_{1}= & \textbf{z}_{0}+\frac{h}{90} K(\frac{\textbf{z}_{0}+\textbf{z}_{1}}{2})(7\nabla H(\textbf{z}_{0})+32\nabla H(\frac{3\textbf{z}_{0}+\textbf{z}_{1}}{4})\\
& +12\nabla H(\frac{\textbf{z}_{0}+\textbf{z}_{1}}{2})+32\nabla H(\frac{\textbf{z}_{0}+3\textbf{z}_{1}}{4})+7\nabla H(\textbf{z}_{1})).\nonumber
\end{align}
The method \eqref{eq:4s5} is exactly energy-conserving if $H$ is a polynomial Hamiltonian of degree $\nu \leq 4$.
It is obvious that the method is symmetric and, therefore, of order 2. By the mono-implicit character of the method, the computational cost is less than high order implicit Runge-Kutta methods which need the values of every stages.

\begin{remark}\label{remark} 
Actually, we can also use other quadrature formulas to calculate the integrand numerically, for instance, the trapezoid rule, the Simpson's rule and so on. We observe that for any non-polynomial Hamiltonian, energy conservation can be practically obtained, by choosing a suitable quadrature formula such that the errors of the Hamiltonian are up to the round-off error~\cite{Refbi}.
In fact, from the numerical experiments we know that the non-polynomial Hamiltonian can be preserved up to the round-off error of $10^{-15}$ if the Boole's rule is used (see Fig. 2(b)).
\end{remark}

\section{Numerical experiments}
\label{sec:5}

In this section, we numerically test the BDLI method \eqref{eq:4s5}.
The fixed point iterative method is used.

\subsection{\textbf{2D dynamics in a static electromagnetic field}}
\label{sec:5s1}

\begin{figure}
  (a)\includegraphics[width=0.45\textwidth]{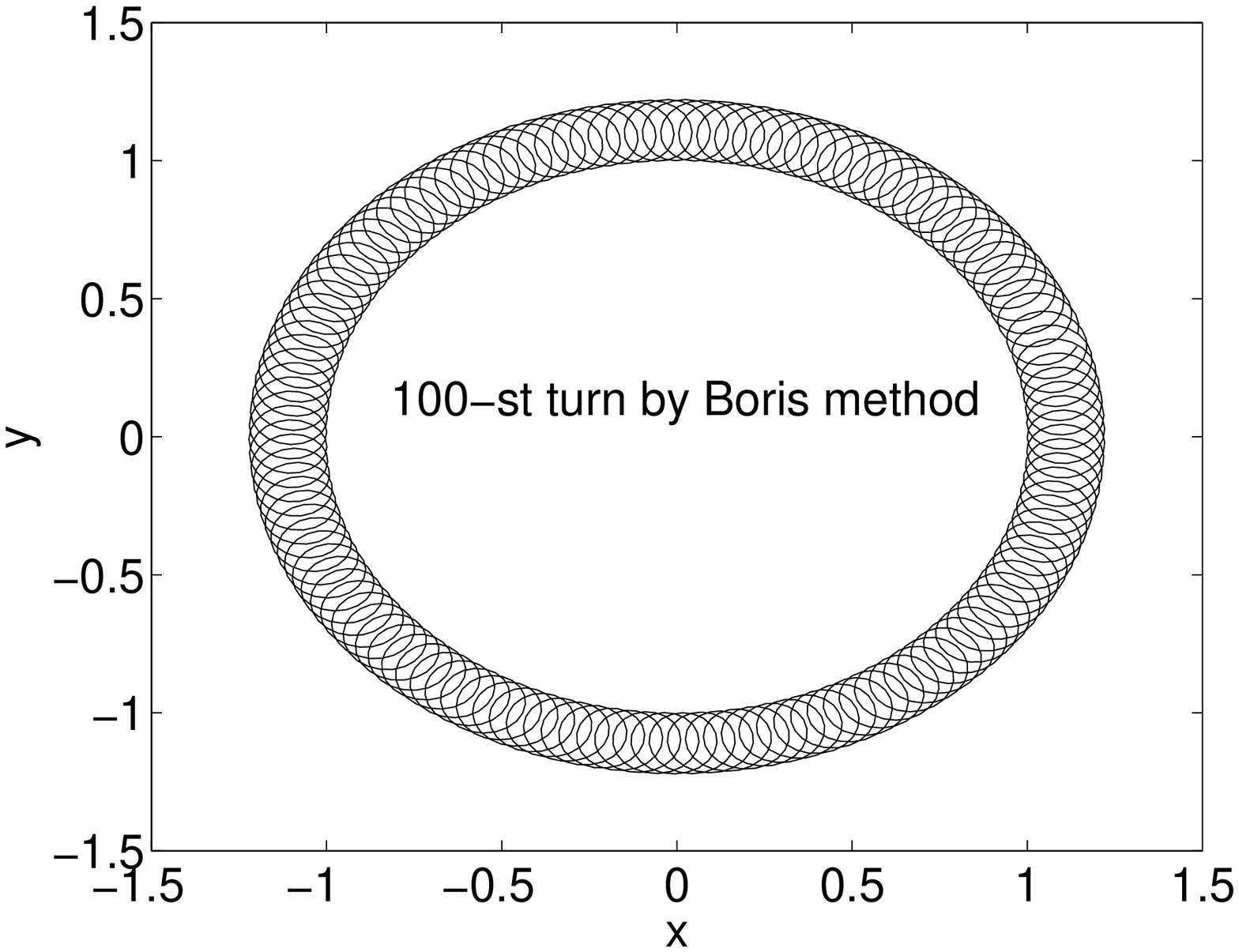}\quad
  (b)\includegraphics[width=0.45\textwidth]{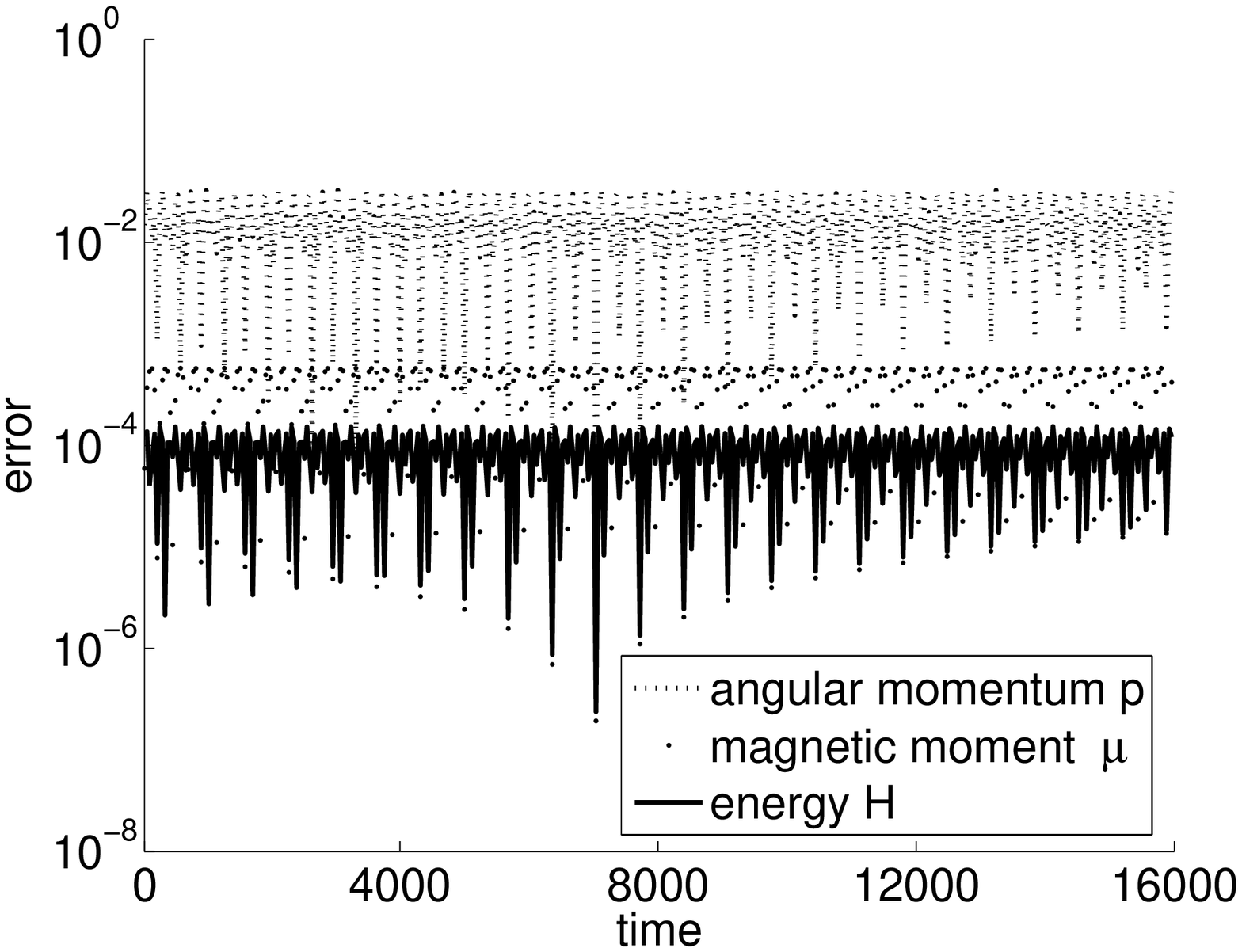}\\
\caption{The Boris method is applied to the simple 2D dynamics with step $h=\pi/10$. (a) The orbit in the 100-st turn; (b) Errors of the angular momentum $p_{\xi}$, the magnetic moment $\mu$ and the energy $H$ for $t\in [0,5\times 10^{4}h]$.}
\label{fig:1}
  (a)\includegraphics[width=0.45\textwidth]{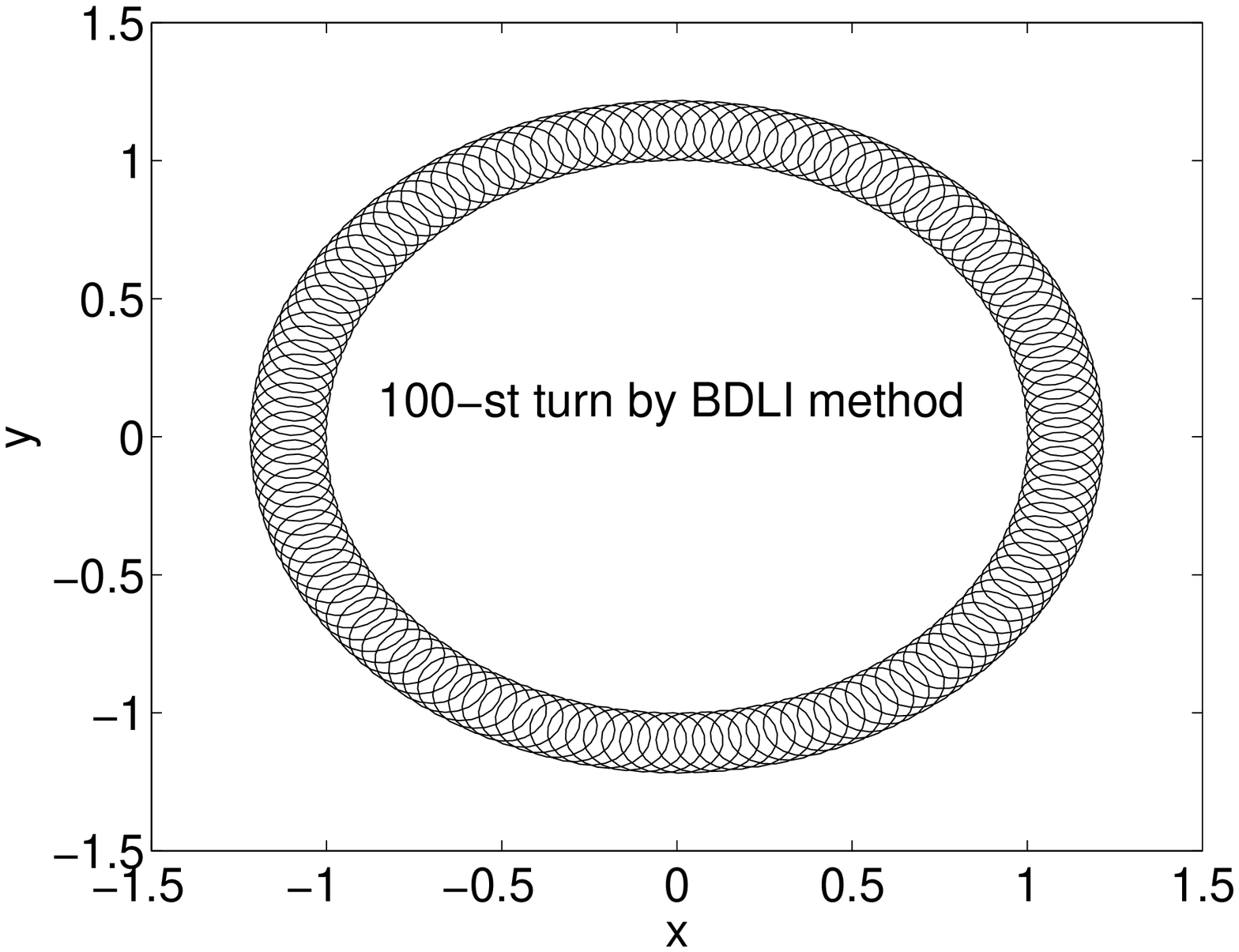}\quad
  (b)\includegraphics[width=0.45\textwidth]{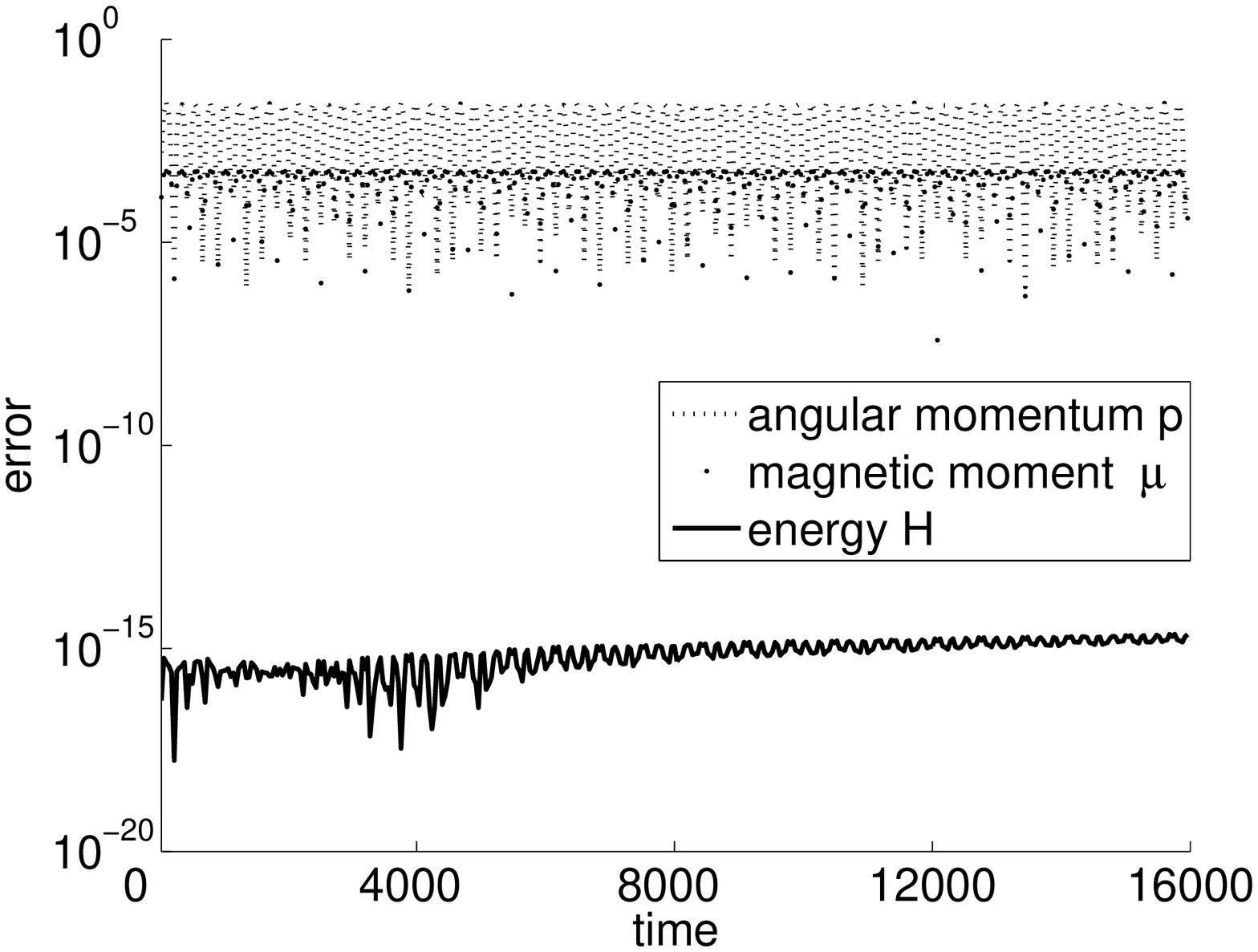}\\
\caption{The BDLI method is applied to the simple 2D dynamics with step $h=\pi/10$. (a) The orbit in the 100-st turn; (b) Errors of the angular momentum $p_{\xi}$, the magnetic moment $\mu$ and the energy $H$ for $t\in [0,5\times 10^{4}h]$.}
\label{fig:2}
\end{figure}

Firstly, we consider a non-polynomial Hamiltonian case, i.e., the 2D dynamics of the charged particle in a static, non-uniform electromagnetic field
\begin{equation}\label{eq:5s1}
\textbf{B}=\nabla\times \textbf{A}=R\textbf{e}_{z},\quad  \textbf{E}=-\nabla \varphi=\frac{10^{-2}}{R^{3}}(x\textbf{e}_{x}+y\textbf{e}_{y}),
\end{equation}
where the potentials are chosen to be $A=\frac{R^{2}}{3}R\textbf{e}_{\xi}$, $\varphi=\frac{10^{-2}}{R}$ in cylindrical coordinates $(R,\xi,z)$ with $R=\sqrt{x^{2}+y^{2}}$. In this example, the physical quantities are normalized by the system size a, the characteristic magnetic field $\textbf{B}_{0}$.
It is known that the energy
\begin{equation}\label{eq:5s2}
H =\frac{1}{2}\textbf{v}\cdot \textbf{v} + \frac{10^{-2}}{R}
\end{equation}
is a constant of motion. As the given electromagnetic field \eqref{eq:5s1} changes slowly with respect to the spatial period of the motion, the magnetic moment $\mu=\frac{\textbf{v}_{\bot}^{2}}{2R}$ is an adiabatic invariant, where $\textbf{v}_{\bot}$ is the component of $\textbf{v}$ perpendicular to $\textbf{B}$.

Starting from the initial position $\textbf{x}_{0}=[0, 0.1, 0]^{T}$ with the initial velocity $\textbf{v}_{0}=[0.1, 0.01, 0]^{T}$, the analytic orbit of the charged particle is a spiraling circle with constant radius. The large circle corresponds to the $\nabla \cdot \textbf{B}$ drift and the $\textbf{E}\times \textbf{B}$ drift of the guiding center, and the small circle is the fast-scale gyro-motion.

We first apply the Boris method and the numerical results are shown in Fig. 1. The step size is $h =\pi/10$ which is the $1/20$ of the characteristic gyro-period.
Fig. 1(a) shows the numerical orbits of the Boris method in the 100-st turn.
In Fig. 1(b), we illustrate the errors of the angular momentum $p_{\xi}$, the magnetic moment $\mu$ and the energy $H$ for $t\in [0,5\times 10^{4}h]$.
It is observed that the errors are bounded over a long integration time.

Next we text the BDLI method \eqref{eq:4s5} with the same step and the numerical results are shown in Fig. 2.
Fig. 2(a) shows the numerical orbits of the BDLI method in the 100-st turn.
Fig. 2(b) shows errors of the invariants for $t\in [0,5\times 10^{4}h]$.
It can be observed that the errors of the angular momentum $p_{\xi}$ and the magnetic moment $\mu$ are also bounded over a long integration time.
What is more, the BDLI method can preserve the non-polynomial Hamiltonian \eqref{eq:5s2} up to the round-off error~\cite{Ref35}.

\subsection{\textbf{2D dynamics in an axisymmetric tokamak geometry}}
\label{sec:5s2}

In this subsection, we consider a polynomial Hamiltonian case, i.e., the motion of a charged particle in a 2-dimensional axisymmetric tokamak geometry without inductive electric field.
The magnetic field in the toroidal coordinates $(r, \theta, \xi)$ is expressed as
\begin{equation}\label{eq:5s3}
\textbf{B}=\frac{B_{0}r}{qR}\textbf{e}_{\theta}+\frac{B_{0}R_{0}}{R}\textbf{e}_{\xi}
\end{equation}
where $B_{0}=1$, $R_{0}=1$, and $q =2$ are constant with their usual meanings. The corresponding vector potential $\textbf{A}$ is chosen to be
\begin{equation}\label{eq:5s4}
\textbf{A}=\frac{z}{2R}\textbf{e}_{R}+\frac{(1-R)^{2}+z^{2}}{4R}\textbf{e}_{\xi}+\frac{\ln R}{2}\textbf{e}_{z}.
\end{equation}
In this example, the physical quantities are normalized by the system size a, characteristic magnetic field $B_{0}$. In the absence of electric field, the energy of the system $H=\frac{1}{2}\textbf{v}\cdot \textbf{v}$ is a polynomial Hamiltonian. With the conserved quantities, the solution orbit projected on
$(R,z)$ space forms a closed orbit.

To apply the discrete line integral method, we transform $\textbf{B}$ \eqref{eq:5s3} into the Cartesian coordinates $(x, y, z)$ which is
\begin{equation}\label{eq:5s5}
\textbf{B}=-\frac{2y+xz}{2R^{2}}\textbf{e}_{x}+\frac{2x-yz}{2R^{2}}\textbf{e}_{y}+\frac{R-1}{2R}\textbf{e}_{y}.
\end{equation}
Starting with the initial position $\textbf{x}_{0}=[1.05, 0, 0]^{T}$ and the initial velocity $\textbf{v}_{0}=[0, 4.816\times 10^{-4}, 2.059\times 10^{-3}]^{T}$, the orbit projected on $(R, z)$ space is a banana orbit, and it will turn into a transit orbit when the initial velocity is changed to $\textbf{v}_{0}=[0, 2\times 4.816\times 10^{-4}, 2.059\times 10^{-3}]^{T}$. Setting the step size $h =\pi/10$ which is the $1/20$ of the gyro-period, we adopt the BDLI method.
The simulation over $5\times 10^{4}$ steps is shown in Fig. 3. Fig. 3(a) shows the banana orbit solution and Fig. 3(a) shows the transit orbit solution.
It is observed that the BDLI method gives correct orbits over a very long integration time. In the case without an electric field, the BDLI method can also preserve the energy H exactly. In Fig. 3(c) we display the errors of the energy. It can be observed that the errors of the energy are up to the round-off error of $10^{-18}$.
Due to the iterative method which is not energy-preserving, a drift is present in the numerical Hamiltonian.
The above experiments show the superiority of the BDLI method for solving the Lorentz force system with a static electromagnetic field.

\begin{figure}
  (a)\includegraphics[width=0.45\textwidth]{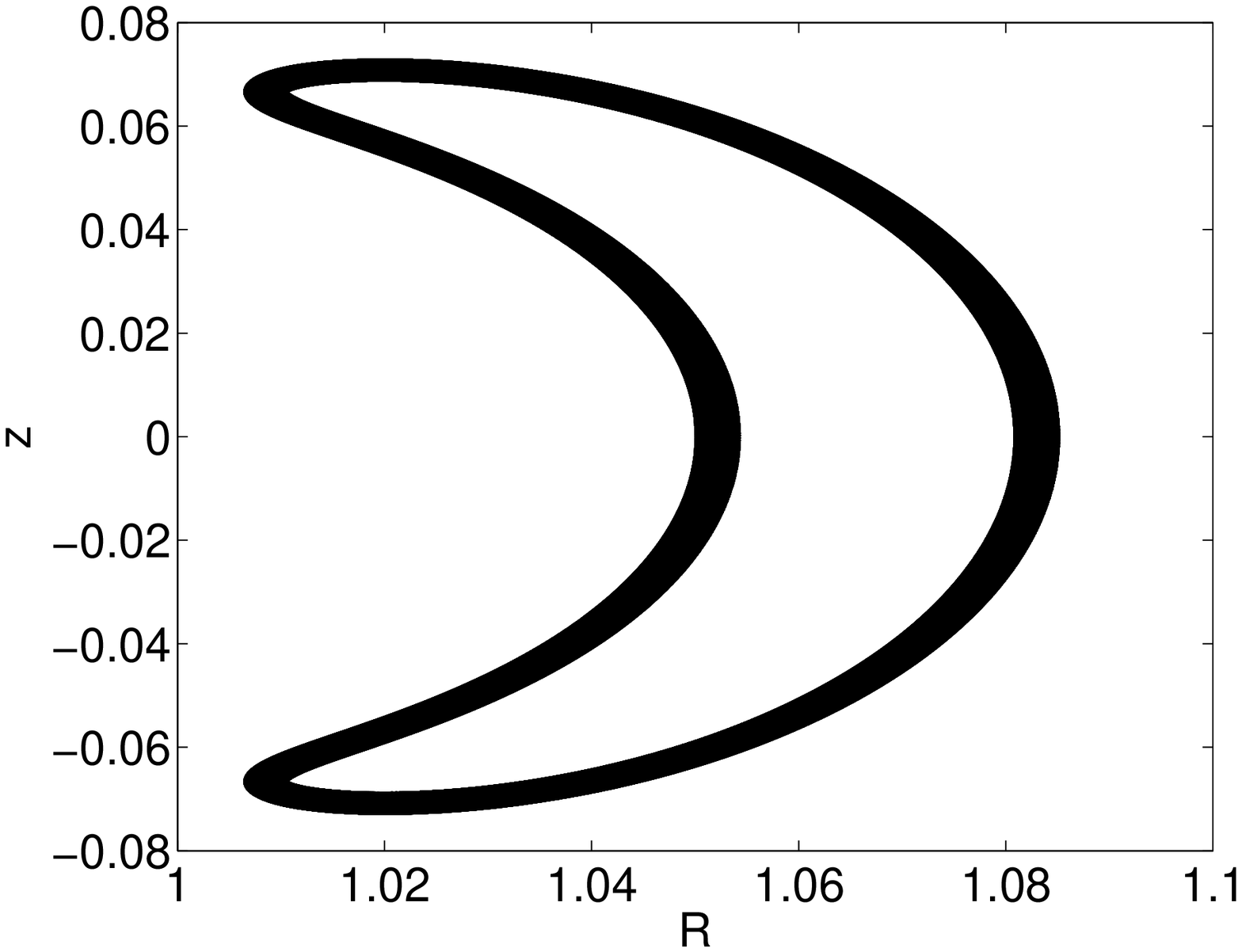}\quad
  (b)\includegraphics[width=0.45\textwidth]{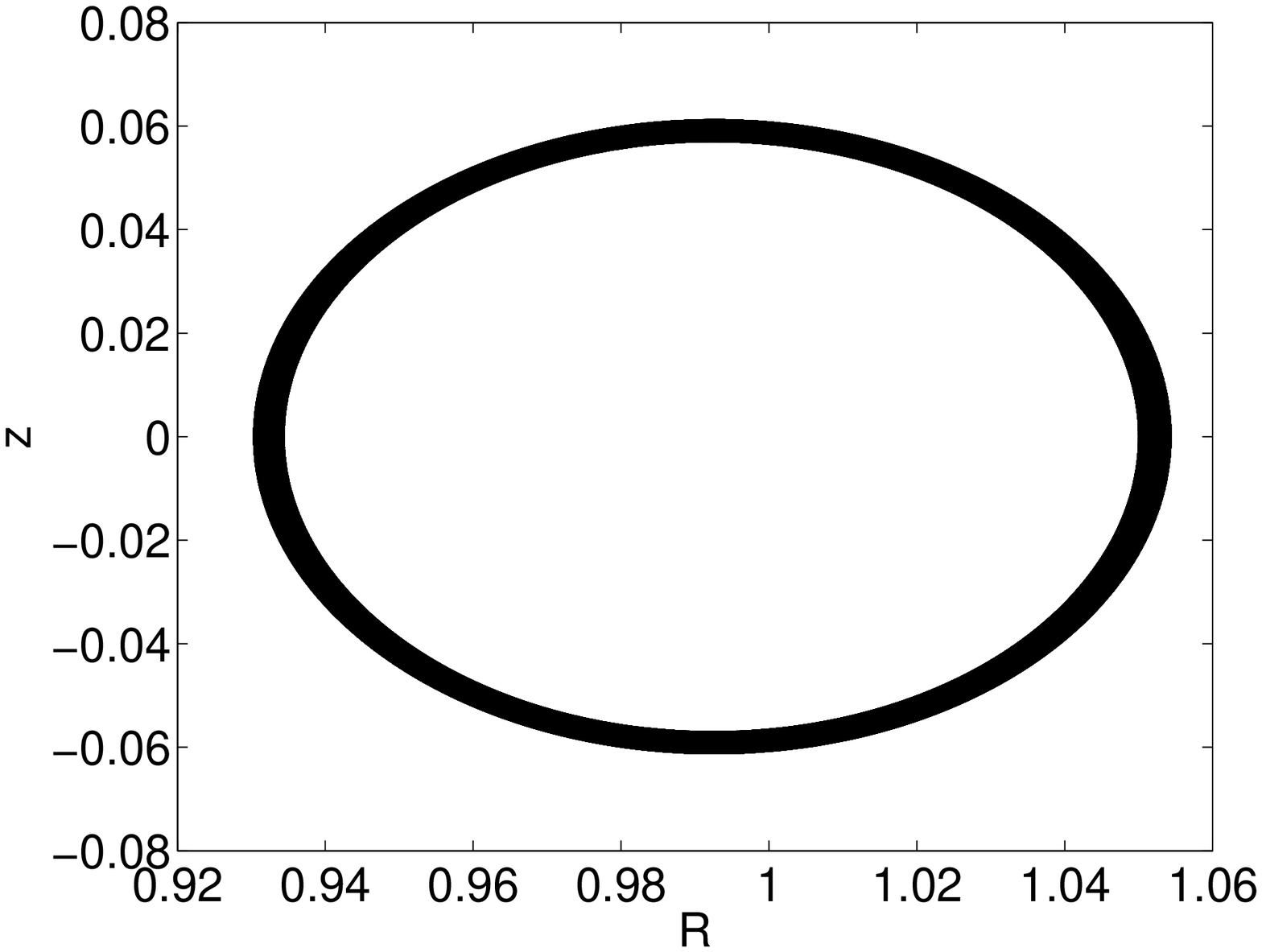}\\
  (c)\includegraphics[width=0.45\textwidth]{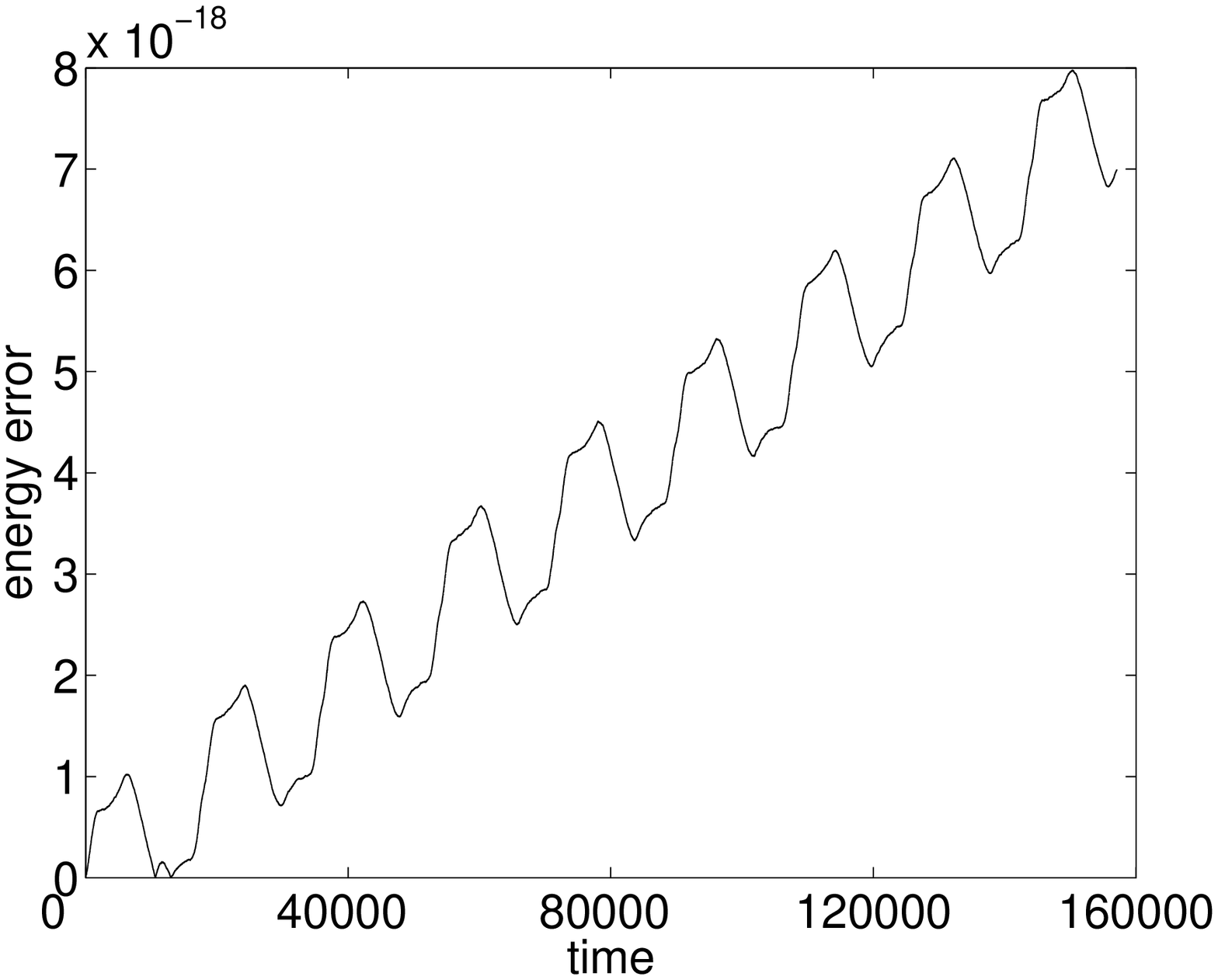}\\
\caption{Numerical solution of the BDLI method with step $h=\pi/10$ for $t\in [0,5\times 10^{5}h]$. (a) Banana orbit; (b) Transit orbit; (c) Energy preservation.}
\label{fig:3}
\end{figure}

\section{Conclusions}
\label{sec:6}

In this paper, we propose the BDLI method for the Lorentz force system which is written as a non-canonical Hamiltonian system in the coordinates $(\textbf{x},\textbf{v})$.
The numerical results show
that, for polynomial Hamiltonian case, the method can conserve the energy exactly, and for non-polynomial Hamiltonian case, the method can also
conserve the energy up to round-off error.
Due to the energy-preserving property,
the BDLI method has good performance over long simulation time.

\vspace{0.3cm}
\hspace{-0.5cm}{\bf Acknowledgements}\\
This work is supported by the Jiangsu Collaborative Innovation Center for Climate
Change, the National Natural Science Foundation of China (Grant Nos. 11271195,
41231173) and the Priority Academic Program Development of Jiangsu Higher Education Institutions.

\vspace{0.5cm}
\hspace{-0.5cm}


{\bf \large Reference}
\begin{thebibliography}{}


\bibitem{Ref1}
R.K. Dodd, J.C. Eibeck, J.D. Gibbon, H.C. Morris, Solitons and Nonlinear Wave Equation, Academic
Press, 1982.
\bibitem{Ref2}
A. Hasegawa, Optical Solitons in Fibers, Springer, Berlin, 1989.
\bibitem{Ref3}
K. Feng, Difference schemes for Hamiltonian formalism and symplectic geometry, J. Comput. Math., 4 (1986) 279-289.
\bibitem{Ref4}
K. Feng, M. Z. Qin, Symplectic Geometric Algorithms for Hamiltonian Systems,
 Springer-Verlag/Zhejiang Science
and Technology Publishing House, Berlin/Hangzhou, 2010.
\bibitem{Refhnw}
E. Hairer, S.P. N$\phi$rsett, G. Wanner, Solving Ordinary Differential Equations. I-Nonstiff Problems, 2nd edition, Springer Ser. Comput. Math., vol. 8, Springer-Verlag, Berlin, 1993.
\bibitem{Refhlw}
E. Hairer, C. Lubich, G. Wanner, Geometric Numerical Integration: Structure-Preserving Algorithms for Ordinary Differential Equations, 2nd edition, Springer-Verlag, Berlin, 2006.
\bibitem{Ref5}
Y.S. Wang, Q.H. Li, Y.Z. Song, Two new simple multi-symplectic schemes for the nonlinear Schr\"{o}dinger equation, Chinese Phys. Lett., 25 (2008) 1538-1540.
\bibitem{Ref6}
J.B. Chen, M.Z. Qin, Multi-symplectic fourier pseudospectral method for the nonlinear Schr\"{o}dinger equation, Electron. Trans. Numer. Anal., 12 (2001) 193-204.
\bibitem{Refkhz}
L.H. Kong, J.L. Hong, J.J. Zang, Splitting multi-symplectic integrators for Maxwell's equation, J. Comput. Phys., 229 (2010) 4259-4278.
\bibitem{Refxsc}
Q. Xu, S.H. Song, Y.M. Chen, A semi-explicit multi-symplectic splitting scheme for a 3-coupled nonlinear Schr\"{o}dinger equation, Comput. Phys. Comm. 185 (2014) 1255-1264.
\bibitem{RefBellan}
P.M. Bellan, Fundamentals of Plasma Physics, 1st edition, Cambridge University Press, 2008.
\bibitem{RefLittlejohn}
R.G. Littlejohn, Hamiltonian formulation of guiding center motion, Phys. Fluids, 24 (1981) 1730-1749.
\bibitem{Refmg}
P.J. Morrison, J.M. Greene, Noncanonical Hamiltonian density formulation of hydrodynamics and ideal magnetohydrodynamics, Phys. Rev. Lett., 45 (1980) 790-794.
\bibitem{RefBoris}
J. Boris, in: Proceedings of the Fourth Conference on Numerical Simulation of Plasmas, Naval Research Laboratory, Washington D.C., 1970, pp. 3-67.
\bibitem{Refscpw}
P.H. Stoltz, J.R. Cary, G. Penn, J. Wurtele, Efficiency of a Boris like integration scheme with spatial stepping, Phys. Rev. Spec. Top., Accel. Beams 5 (2002) 094001.
\bibitem{Refbl}
C.K. Birdsall, A.B. Langdon, Plasma Physics via Computer Simulation, Ser. Plasma Phys., Taylor $\&$ Francis, 2005.
\bibitem{Refqzxlst}
H. Qin, S.X. Zhang, J.Y. Xiao, J. Liu, Y.J. Sun, W.M. Tang, Why is Boris algorithm so good?, Phys. Plasmas, 20 (2013) 084503.
\bibitem{Refhslq}
Y. He, Y.J. Sun, J. Liu, H. Qin, Volume-preserving algorithms for charged particle dynamics, J. Comp. Phys., 281 (2015) 135-147.
\bibitem{Refqgt}
H. Qin, X. Guan, W.M. Tang, Variational symplectic algorithm for guiding center dynamics and its application in tokamak geometry, Phys. Plasmas, 16 (2009) 042510.
\bibitem{Refqg}
H. Qin, X. Guan, Variational symplectic integrator for long-time simulations of the guiding-center motion of charged particles in general magnetic fields, Phys. Rev. Lett., 100 (2008) 035006.
\bibitem{RefLittlejohn2}
R.G. Littlejohn, Variational principles of guiding centre motion, J. Plasma Phys., 29 (1983) 111-125.
\bibitem{Ref31}
O. Gonzalez, Time integration and discrete Hamiltonian systems, J. Nonlinear Sci., 6 (1996) 449-467.
\bibitem{Ref32}
T. Matsuo, High-order schemes for conservative or dissipative systems, J. Comput. Appl. Math., 152 (2003) 305-317.
\bibitem{Ref37}
G. R. W. Quispel, D. I. McLaren, A new class of energy-preserving numerical integration methods,
J. Phys. A: Math. Theor.,  41 (2008) 045206.
\bibitem{Ref30}
Y.Z. Gong, J.X, Cai, Y.S. Wang, Some new structure-preserving algorithms for general
multi-symplectic formulations of Hamiltonian PDEs, J. Comp. Phys., 279 (2014) 80-102.
\bibitem{Ref36}
L. Brugnano, F. Iavernaro, D. Trigiante. A note on the efficient implementation of
Hamiltonian BVMs, J. Comput. Appl. Math., 236 (2011) 375-383.
\bibitem{Refip1}
F. Iavernaro, B. Pace. S-Stage trapezoidal methods for the conservation of Hamiltonian
Functions of Polynomial Type, AIP Conf. Proc., 936 (2007) 603-606.
\bibitem{Refbi}
L. Brugnano, F. Iavernaro, Line integral methods
and their application to the numerical solution of
conservative problems, (2013) arXiv:1301.2367v1 [math.NA].
\bibitem{Ref33}
L. Brugnano, F. Iavernaro,  D. Trigiante,  Hamiltonian boundary value methods(Energy preserving discrete line integral
methods), J. Numer. Anal. Ind. Appl. Math., 5 (2010) 17-37.
\bibitem{Ref34}
F. Iavernaro, D. Trigiante, High-order symmetric schemes for the energy conservation of
polynomial Hamiltonian problems, J. Numer. Anal. Ind. Appl. Math., 4 (2009) 87-101.
\bibitem{Ref35}
L. Brugnano, F. Iavernaro, D. Trigiante. A simple framework for the derivation and analysis of effective one-step methods for ODEs, Appl. Math. Comput.,
218 (2012) 8475-8485.


\end{thebibliography}
\end{document}